\documentclass[11pt,reqno,twoside]{amsart}
\usepackage{amsmath}
\usepackage{amssymb,amsfonts,amsthm,indentfirst}

\newtheorem{theo}{Theorem}[section]

\newtheorem{prob}{Problem}[section]
\newtheorem{lem}{Lemma}[section]
\newtheorem{rem}{Remark}[section]
\newtheorem{cor}{Corollary}[section]

\newcommand{\be}{\begin{equation}}
\newcommand{\ee}{\end{equation}}
\newcommand{\bea}{\begin{eqnarray}}
\newcommand{\eea}{\end{eqnarray}}
\newcommand{\beas}{\begin{eqnarray*}}
\newcommand{\eeas}{\end{eqnarray*}}

\begin{document}
\title[Pseudo-parallel Normal Jacobi Operator...]{Real Hypersurfaces with
Pseudo-parallel Normal Jacobi Operator in Complex Two-Plane Grassmannians}
\thanks{This work was supported in part by the UMRG research grant (Grant
No. RG163-11AFR)}
\address{De, A.\\
Institute of mathematical Sciences, University of Malaya, 50603 Kuala
Lumpur, Malaysia}
\email{de.math@gmail.com}
\author{Avik De}
\author{Tee-How Loo}
\address{Loo, T. H. \\
Institute of Mathematical Sciences, University of Malaya, 50603 Kuala
Lumpur, Malaysia.}
\email{looth@um.edu.my}

\begin{abstract}
The objective of the present paper is to prove the non-existence of real
hypersurface with pseudo-parallel normal Jacobi operator in complex
two-plane Grassmannians. As a corollary, we show that there does not exist
any real hypersurface with semi-parallel or recurrent normal Jacobi operator
in complex two-plane Grassmannians. This answers a question considered in
[Monatsh Math, 172 (2013), 167-178] in negative.
\end{abstract}

\maketitle

\setcounter{page}{1} \setlength{\unitlength}{1mm}\baselineskip .45cm %
\pagenumbering{arabic} \numberwithin{equation}{section}




\section{Introduction.}

A complex two-plane Grassmannian $G_{2}(\mathbb{C}_{m+2})$ is the set of all
complex two-dimensional linear subspaces of $\mathbb{C}_{m+2}$. It is the
unique compact, irreducible Riemannian symmetric space with positive scalar curvature, equipped with both a Kaehler structure $J$ and a quaternionic Kaehler structure $\mathfrak{J}$
not containing $J$ \cite{berndt}.

Typical examples of real hypersurfaces $M$ in $G_2(\mathbb{C}_{m+2})$ are
tubes around $G_2(\mathbb{C}_{m+1})$ and tubes around $\mathbb{H}P^{n}$ in $%
G_2(\mathbb{C}_{2n+2})$ . These two classes of real hypersurfaces possess a
number of interesting geometric properties. Characterizing them or
subclasses of them under certain nice geometric conditions has been one of
the main focus of researchers in the theory of real hypersurfaces in $G_2(%
\mathbb{C}_{m+2})$. One of the foremost results along this line was obtained by Berndt
and Suh \cite{berndt-suh}, they characterized these two classes of real
hypersurfaces under the invariance of vector bundles $JT^\perp M$ and $%
\mathfrak{J}T^\perp M$ over the real hypersurfaces $M$ under the shape
operator $A$ of $M$, where $T^\perp M$ is the normal bundle of $M$.

On the other hand, these structures $J$ and $\mathfrak{J}$ of $G_2(\mathbb{C}%
_{m+2})$ significantly impose restrictions on the geometry of its real
hypersurfaces. For instance, there does not exist any parallel real
hypersurface \cite{suh6}. Determining the existence (or non-existence) of
real hypersurfaces in $G_2(\mathbb{C}_{m+2})$ satisfying certain geometric
properties has also became another main research topic in this theory. The main
objective of this paper is to prove the non-existence of real hypersurfaces
in $G_2(\mathbb{C}_{m+2})$ with pseudo-parallel normal Jacobi operator.

Recall that the \emph{normal Jacobi operator} $\hat R_N$, for a hypersurface 
$M$ in a Riemannian manifold, is defined as $\hat{R}_{N}(X)=\hat{R}(X,N)N$,
for any vector $X$ tangent to $M$, where $\hat{R}$ is the curvature tensor
of the ambient space and $N$ is a unit vector normal to $M$ \cite{ber}.

Let $M$ be an orientable real hypersurface isometrically immersed in $G_2(%
\mathbb{C}_{m+2})$. Denote by $(\phi,\xi,\eta)$ the almost contact structure
on $M$ induced by $J$, $(\phi_a,\xi_a,\eta_a)$, $a\in\{1,2,3\}$, the local
almost $3$-structure on $M$ induced by $\mathfrak{J}$ and $\mathfrak{D}%
^\perp= \mathfrak{J}T^\perp M$. The real hypersurface $M$ is said to be 
\emph{Hopf} if $AJT^\perp M\subset JT^\perp M$, or equivalently, the Reeb
vector field $\xi$ is principal with principal curvature $\alpha$.

P\'{e}rez et al. \cite{perez}\ studied the real hypersurfaces in $G_{2}(%
\mathbb{C}_{m+2})$, in which the normal Jacobi operator commutes with both
the shape operator and the structure tensor $\phi $. In \cite{jeong3} Jeong
et al. proved the following:
\begin{theo}[\protect\cite{jeong3}]
\label{thm:suh} There does not exist any connected Hopf hypersurface in
complex two-plane Grassmannians $G_{2}(\mathbb{C}_{m+2}),m\geq 3$, with
parallel normal Jacobi operator.
\end{theo}
Machado et al. \cite{machado} proved the non-existence of Hopf hypersurfaces
in $G_{2}(\mathbb{C}_{m+2})$ with Codazzi type $\hat{R}_{N}$ under certain
conditions on the $\mathfrak{D}$- and $\mathfrak{D}^{\bot }$-component of $%
\xi $ . Later, the non-existence of Hopf hypersurfaces in $G_{2}(\mathbb{C}%
_{m+2})$ whose normal Jacobi operator is $(\mathbb{R}\xi\cup \mathfrak{D}%
^{\bot })$-parallel was proved \cite{jeong}. In \cite{jeong 2}, Suh and
Jeong investigated real hypersurfaces in $G_{2}(\mathbb{C}_{m+2})$ with $%
L_{\xi }\hat R_{N}=0$, and proved the non-existence of such real
hypersurfaces under the condition either $\xi \in \mathfrak{D}$ or $\xi \in 
\mathfrak{D}^{\bot }$. They also proved the non-existence of Hopf
hypersurfaces with Lie parallel normal Jacobi operator, that is, $L_{X}\hat
R_{N}=0$ in $G_{2}(\mathbb{C}_{m+2})$ \cite{jeong1}.

A real hypersurface $M$ is said to have \emph{recurrent} normal Jacobi
operator if $(\hat\nabla_X\hat R_N)Y=\omega(X)\hat R_NY$, for some 1-form $%
\omega$. In \cite{suh rec}, Jeong et al. generalized Theorem \ref{thm:suh}
and proved the following:
\begin{theo}[\protect\cite{suh rec}]
\label{recc} There does not exist any connected Hopf hypersurface in complex
two-plane Grassmannians $G_{2}(\mathbb{C}_{m+2})$, with recurrent normal
Jacobi operator.
\end{theo}
Deprez \cite{deprez2} first studied a submanifold $M$ in a Riemannian
manifold whose second fundamental form $h$ satisfies $\bar{R}\cdot h=0$,
where $\bar{R}$ is the curvature tensor corresponding to the van der
Waerden-Bortolotti connection. Such submanifolds are said to be \emph{%
semi-parallel}. In \cite{ortega}, Ortega proved that there does not exist
any semi-parallel real hypersurface in a non-flat complex space form.

A $(1,p)$-tensor $T$, in a Riemannian manifold $M$ with Riemannian curvature
tensor ${R}$ is said to be \emph{pseudo-parallel}, if it satisfies ${R}(X,Y)
T=f\{(X\wedge Y)T\}$, for some function $f$, where 
\begin{align*}
(X\wedge Y)Z:=\langle Y,Z\rangle X-\langle X,Z\rangle Y,
\end{align*}
and 
\begin{align*}
\{(X\wedge Y)T\}(X_1,\cdots,X_p):=(X\wedge Y)T(X_1,\cdots,X_p) \\
-\sum_{j=1}^pT(X_1,\cdots (X\wedge Y)X_j,\cdots, X_p).
\end{align*}
for any $X,Y,Z,X_1,\cdots,X_p\in TM$.

The notion of pseudo-parallel submanifolds, that is, submanifolds with
pseudo-parallel second fundamental forms, can be considered as a
generalization of semi-parallel submanifolds. Asperti et al. \cite{asperti}
classified all pseudo-parallel hypersurfaces in space forms as quasi-umbilic
hypersurfaces or cyclids of duplin. The classification of pseudo-parallel real hypersurfaces in a non-flat complex space form was obtained in
\cite{lobos}.

Recently, Panagiotidou and Tripathi \cite{trip} studied Hopf hypersurfaces
with semi-parallel normal Jacobi operator in $G_{2}(\mathbb{C}_{m+2})$ and
proved the following:
\begin{theo}[\protect\cite{trip}]
$\label{tripa}$ There does not exist any connected Hopf hypersurface $M$ in $%
G_{2}(\mathbb{C}_{m+2}),m\geq 3$,\ equipped with semi-parallel normal Jacobi
operator, if $\alpha\neq0$ and $\mathfrak{D}$- or $\mathfrak{D}^{\bot }$%
-component of the Reeb vector field $\xi$ is invariant by the shape operator 
$A$.
\end{theo}
One of the challenges in the theory of real hypersurfaces $M$ in $G_{2}(\mathbb{C%
}_{m+2})$ is handling lengthy and complicated expressions resulting from the complexity of the geometric structures on $M$, induced by the Kaehler and the quaternionic Kaehler structure of $G_{2}(\mathbb{C}_{m+2})$. 
For technical reasons, certain additional restrictions like 
$M$ being Hopf, having non-vanishing geodesic Reeb flow etc have often been imposed while dealing with real hypersurfaces in $G_{2}(\mathbb{C}_{m+2})$. It would be interesting to see if any nice results on real hypersurfaces of $G_{2}(\mathbb{C}_{m+2})$ can be obtained without these restrictions. 

Motivated by Theorem~\ref{thm:suh}, Theorem \ref{recc} and Theorem~%
\ref{tripa}, a question arises naturally:

\begin{prob}\label{bad}
Does there exist real hypersurface in $G_{2}(\mathbb{C}_{m+2})$ with
parallel, recurrent, semi-parallel or pseudo-parallel normal Jacobi operator?
\end{prob}

We first prove the following:
\begin{theo}
$\label{mainthm}$ There does not exist any real hypersurface with
pseudo-parallel normal Jacobi operator in $G_{2}(\mathbb{C}_{m+2})$, $m\geq3$.
\end{theo}
\begin{rem}
It is worthwhile to note that no additional condition has been imposed in  the above theorem.
\end{rem}

Since a semi-parallel tensor is always pseudo-parallel with the associated
function $f=0$, we state that

\begin{cor}
\label{cor1.5} There does not exist any real hypersurface in $G_{2}(\mathbb{C%
}_{m+2})$, $m\geq3$, equipped with semi-parallel normal Jacobi operator.
\end{cor}
In the last section we prove the non-existence of real hypersurfaces with recurrent or parallel normal Jacobi operator in $G_{2}(\mathbb{C}_{m+2})$, $m\geq3$ (see Corollary \ref{4.1}). Thus Problem \ref{bad} has been solved completely.


\section{Real hypersurfaces in $G_2(\mathbb{C}_{m+2})$}

In this section we state some structural equations as well as some known
results in the theory of real hypersurfaces in $G_2(\mathbb{C}_{m+2})$. A
thorough study on the Riemannian geometry on $G_{2}(\mathbb{C}_{m+2})$ can
be found in \cite{berndt}. Denote by $\langle ,\rangle $ the Riemannian
metric, $J$ the Kaehler structure and $\mathfrak{J}$ the quaternionic
Kaehler structure on $G_{2}(\mathbb{C}_{m+2})$. For each $x\in G_{2}(\mathbb{%
C}_{m+2})$, we denote by $\{J_{1},J_{2},J_{3}\}$ a canonical local basis of $%
\mathfrak{J}$ on a neighborhood $\mathcal{U}$ of $x$ in $G_{2}(\mathbb{C}%
_{m+2})$, that is, each $J_{a}$ is a local almost Hermitian structure such
that 
\begin{equation}
J_{a}J_{a+1}=J_{a+2}=-J_{a+1}J_{a},\quad a\in \{1,2,3\}.\label{eqn:quarternion}
\end{equation}
Here, the index is taken modulo three. Denote by $\hat{\nabla}$ the
Levi-Civita connection of $G_{2}(\mathbb{C}_{m+2})$. There exist local $1$%
-forms $q_{1}$, $q_{2}$ and $q_{3}$ such that 
\begin{equation*}
\hat{\nabla}_{X}J_{a}=q_{a+2}(X)J_{a+1}-q_{a+1}(X)J_{a+2}
\end{equation*}%
for any $X\in T_{x}G_{2}(\mathbb{C}_{m+2})$, that is, $\mathfrak{J}$ is
parallel with respect to $\hat{\nabla}$. The Kaehler structure $J$ and
quaternionic Kaehler structure $\mathfrak{J}$ is related by 
\begin{equation}
JJ_{a}=J_{a}J;\quad \text{Trace}{(JJ_{a})}=0,\quad a\in \{1,2,3\}.\label{eqn:JJa}
\end{equation}%

The Riemannian curvature tensor $\hat R$ of $G_2(\mathbb{C}_{m+2})$ is
locally given by 
\begin{align}  \label{eqn:hatR}
\hat R(X,Y)Z=&\langle Y,Z\rangle X-\langle X,Z\rangle Y+\langle JY,Z\rangle
JX-\langle JX,Z\rangle JY-2\langle JX,Y\rangle JZ  \notag \\
&+\sum_{a=1}^3\{\langle J_aY,Z\rangle J_aX-\langle J_aX,Z\rangle
J_aY-2\langle J_aX,Y\rangle J_aZ  \notag \\
&+\langle JJ_aY,Z\rangle JJ_aX-\langle JJ_aX,Z\rangle JJ_aY\}.
\end{align}
for all $X$, $Y$ and $Z\in T_xG_2(\mathbb{C}_{m+2})$.

For a nonzero vector $X\in T_xG_2(\mathbb{C}_{m+2})$, we denote by $\mathbb{C%
}X=\text{Span}\{X,JX\}$, $\mathfrak{J}X=\{J^{\prime }X|J^{\prime }\in%
\mathfrak{J}_x\}$, $\mathbb{H}X=\mathbb RX\oplus\mathfrak{J}X$, and $%
\mathbb{H}\mathbb{C}X$ the subspace spanned by $\mathbb{H}X$ and $\mathbb{H}%
JX$. If $JX\in\mathfrak{J}X$, we denote by $\mathbb{C}^\perp X$ the
orthogonal complement of $\mathbb{C}X$ in $\mathbb{H}X$.

Let $M$ be an oriented real hypersurface isometrically immersed in $G_2(%
\mathbb{C}_{m+2})$, $m\geq3$, $N$ a unit normal vector field on $M$. The
Riemannian metric on $M$ is denoted by the same $\langle,\rangle$. A
canonical local basis $\{J_1,J_2,J_3\}$ of $\mathfrak{J}$ on $G_2(\mathbb{C}%
_{m+2})$ induces a local almost contact metric $3$-structure $%
(\phi_a,\xi_a,\eta_a,\langle,\rangle)$ on $M$ by 
\begin{align*}
J_aX=\phi_a X+\eta_a(X)N, \quad J_aN=-\xi_a, \quad
\eta_a(X)=\langle\xi_a,X\rangle
\end{align*}
for any $X\in TM$. It follows from (\ref{eqn:quarternion}) that 
\begin{align*}
&\phi_a\phi_{a+1}-\xi_a\otimes\eta_{a+1}=\phi_{a+2}=-\phi_{a+1}\phi_a+%
\xi_{a+1}\otimes\eta_a \\
&\phi_a\xi_{a+1}=\xi_{a+2}=-\phi_{a+1}\xi_a.
\end{align*}
Denote by $(\phi, \xi,\eta,\langle,\rangle)$ the almost contact metric
structure on $M$ induced by $J$, that is, 
\begin{align*}
JX=\phi X+\eta(X)N, \quad JN=-\xi, \quad \eta(X)=\langle\xi,X\rangle.
\end{align*}
The vector field $\xi$ is known as the \emph{Reeb vector field}.

It follows from (\ref{eqn:JJa}) that the two structures $(\phi,\xi,\eta,%
\langle,\rangle)$ and $(\phi_a,\xi_a,\eta_a,\langle,\rangle)$ can be related
as follows 
\begin{align*}
\phi_a\phi-\xi_a\otimes\eta=\phi\phi_a-\xi\otimes\eta_a; \quad
\phi\xi_a=\phi_a\xi.
\end{align*}
Denote by $\nabla$ the Levi-Civita connection and $A$ the shape operator on $%
M$. Then 
\begin{align*}
&(\nabla_{X} \phi)Y=\eta(Y)AX-\langle AX,Y\rangle\xi, \quad \nabla_X \xi =
\phi AX  \label{eqn:delxi} \\
&(\nabla_{X}\phi_a)Y=\eta_a(Y)AX-\langle
AX,Y\rangle\xi_{a}+q_{a+2}(X)\phi_{a+1}Y-q_{a+1}(X)\phi_{a+2}Y \\
&\nabla_X \xi_a = \phi_a AX+q_{a+2}(X)\xi_{a+1}-q_{a+1}(X)\xi_{a+2}
\end{align*}
for any $X,Y\in TM$.

Corresponding to each canonical local basis $\{J_1,J_2,J_3\}$ of $\mathfrak{J%
}$, we define a local endomorphism $\theta_a$ on $TM$ by 
\begin{align*}
\theta_aX:=\tan(J_aJX)=\phi_a\phi X-\eta(X)\xi_a.
\end{align*}
Some properties of $\theta_a$ are given in the following:

\begin{lem}[\protect\cite{loo}]
\label{lem:theta}

\begin{enumerate}
\item[(a)] $\theta_a$ is symmetric,

\item[(b)] $\text{Trace}\,{\theta_a}=\eta(\xi_a)$,

\item[(c)] $\theta_a^2X=X+\eta_a(X)\phi\xi_a$, for all $X\in TM$,

\item[(d)] $\theta_a\xi=-\xi_a; \quad \theta_a\xi_a=-\xi; \quad
\theta_a\phi\xi_a=\eta(\xi_a)\phi\xi_a$,

\item[(e)] $\theta_a\xi_{a+1}= \phi\xi_{a+2}=-\theta_{a+1}\xi_a$,

\item[(f)] $\theta_a\phi\xi_{a+1}=-\xi_{a+2}+\eta(\xi_{a+1})\phi\xi_a$,

\item[(g)] $\theta_{a+1}\phi\xi_a= \xi_{a+2}+\eta(\xi_a)\phi\xi_{a+1}$.
\end{enumerate}
\end{lem}

For each $x\in M$, we define a subspace $\mathcal{H}^\perp$ of $T_xM$ by 
\begin{equation*}
\mathcal{H}^\perp:=\text{Span}\{\xi,\xi_1,\xi_2,\xi_3,\phi\xi_1,\phi\xi_2,%
\phi\xi_3\}.
\end{equation*}
Let $\mathcal{H}$ be the orthogonal complement of $\mathbb{H}\mathbb{C}\xi$
in $T_xG_2(\mathbb{C}_{m+2})$. Then $T_xM=\mathcal{H}\oplus\mathcal{H}^\perp$%
.

\begin{lem}[\protect\cite{loo}]
\label{lem:theta_2} Let $\mathcal{H}_a(\varepsilon)$ be the eigenspace
corresponding to eigenvalue $\varepsilon$ of $\theta_a$. Then 

\begin{enumerate}
\item[(a)] ${\theta_a}_{|\mathcal{H}}$ has two eigenvalues $\varepsilon=\pm1$%
,

\item[(b)] $\phi\mathcal{H}_a(\varepsilon)=\phi_a\mathcal{H}_a(\varepsilon)=%
\mathcal{H}_a(\varepsilon)$,

\item[(c)] $\theta_b\mathcal{H}_a(\varepsilon)=\mathcal{H}_a(-\varepsilon)$,
for $a\neq b$.

\item[(d)] $\phi_b\mathcal{H}_a(\varepsilon)=\mathcal{H}_a(-\varepsilon)$,
for $a\neq b$.
\end{enumerate}
\end{lem}

Let $R$ be the curvature tensor of $M$. It follows from (\ref{eqn:hatR})
that the equation of Gauss is given by 
\begin{align}
R(X,Y)Z &=\langle Y,Z\rangle X-\langle X,Z\rangle Y+\langle \phi Y,Z\rangle
\phi X-\langle \phi X,Z\rangle \phi Y-2\langle \phi X,Y\rangle \phi Z  \notag
\\
&+\sum_{a=1}^{3}\{\langle \phi _{a}Y,Z\rangle \phi _{a}X-\langle \phi
_{a}X,Z\rangle \phi _{a}Y-2\langle \phi _{a}X,Y\rangle \phi _{a}Z  \notag \\
&+\langle \theta _{a}Y,Z\rangle \theta _{a}X-\langle \theta _{a}X,Z\rangle
\theta _{a}Y\}+\langle AY,Z\rangle AX-\langle AX,Z\rangle AY,  \label{ap1}
\end{align}
for any $X,Y,Z\in TM$.

The normal Jacobi operator $\hat{R}_{N}$ is given by 
\begin{equation}
\hat{R}_{N}X=X+3\eta (X)\xi +\sum_{a=1}^{3}\{3\eta _{a}(X)\xi _{a}-\eta (\xi
_{a})\theta _{a}+\eta _{a}(\phi X)\phi \xi _{a}\}.  \label{eqn:jacobi}
\end{equation}%
for any $X\in TM$.


Finally, we state two known results which we use in the next section:

\begin{lem}[\protect\cite{loo}]
\label{lem:Xua} Let $M$ be a real hypersurface in $G_2(\mathbb{C}_{m+2})$.
If $\xi$ is tangent to $\mathfrak{D}$ then $A\phi\xi_a=0$, for $%
a\in\{1,2,3\} $.
\end{lem}

\begin{lem}[\protect\cite{loo}]
\label{lem:dga} Let $M$ be a real hypersurface in $G_2(\mathbb{C}_{m+2})$.
If $\xi$ is tangent to $\mathfrak{D}$ then $\xi$, $\xi_1,\xi_2,\xi_3,\phi%
\xi_1,\phi\xi_2,\phi\xi_3$ are orthonormal.
\end{lem}


\section{Proof of Theorem~\protect\ref{mainthm}}

Let the normal Jacobi operator $\hat{R}_{N}$ for a real hypersurface $M$ in $%
G_{2}(\mathbb{C}_{m+2})$, be pseudo-parallel. Then we have 
\begin{equation}
\langle (R(X,Y)\hat{R}_{N})Z,W\rangle =f\langle \lbrack (X\wedge Y)\hat{R}%
_{N}]Z,W\rangle ,
\end{equation}
for any $X,Y,Z,W\in TM$, where $f$ is a real-valued function on $M$, This
implies that 
\begin{eqnarray}
&\langle R(X,Y)\hat{R}_{N}Z,W\rangle -\langle \hat{R}_{N}R(X,Y)Z,W\rangle
=f\{\langle Y,\hat{R}_{N}Z\rangle \langle X,W\rangle  \notag \\
&-\langle X,\hat{R}_{N}Z\rangle \langle Y,W\rangle -\langle Y,Z\rangle
\langle \hat{R}_{N}X,W\rangle+\langle X,Z\rangle \langle \hat{R}%
_{N}Y,W\rangle \}.  \label{ap5}
\end{eqnarray}
We consider two cases: $\xi\notin\mathfrak{D}$ at a point $x\in M$; and $%
\xi\in\mathfrak{D}$ on $M$.

\textbf{Case 1: }$\xi\notin\mathfrak{D}$ at a point $x\in M$.

Without loss of generality, we assume $0<\eta(\xi_1)\leq 1$, $%
\eta(\xi_2)=\eta(\xi_3)=0$. Let $\beta$, $\mu\in\mathbb{R}$ and $U\in%
\mathcal{H}_1(1)$, $V\in\mathcal{H}_1(-1)$ be unit vectors such that the $%
\mathcal{H}_1(1)$-component $(A\xi)^+$ and $\mathcal{H}_1(-1)$-component $%
(A\xi)^-$ of $A\xi$ are given by 
\begin{equation*}
(A\xi)^+=\beta U, \qquad (A\xi)^-=\mu V.
\end{equation*}

Since, $\eta(\xi_2)=0=\eta(\xi_3)$, for $Z\in \mathcal{H}_{1}(1)$ and $W\in 
\mathcal{H}_{1}(-1)$, from (\ref{eqn:jacobi}) we have 
\beas
\hat R_NZ=(1-\eta(\xi_1))Z, \quad \hat R_NW=(1+\eta(\xi_1))W.
\eeas
Since $\eta (\xi _{1})\neq 0$, by putting $Z\in \mathcal{H}_{1}(1)$ and $W\in \mathcal{H}_{1}(-1)$ in (\ref{ap5}), we obtain 
\begin{equation}
\langle R(X,Y)Z,W\rangle =f\{\langle Y,Z\rangle\langle X,W\rangle-\langle
X,Z\rangle\langle Y,W\rangle\},  \label{ap}
\end{equation}
for any $X$, $Y\in T_{x}M$, $Z\in \mathcal{H}_{1}(1)$ and $W\in \mathcal{H}%
_{1}(-1)$. In particular, for $X=\xi $ and $Y\perp \xi $, using the Gauss
equation (\ref{ap1}), the above equation becomes 
\begin{eqnarray*}
&-2\sum_{a=1}^{3}\langle \phi \xi _{a},Y\rangle \langle
\phi_aZ,W\rangle+\sum_{a=1}^{3}\{\langle \theta_aY,Z\rangle\langle
\theta_a\xi, W\rangle-\langle \theta_aY,W\rangle\langle
\theta_a\xi,Z\rangle\}  \notag \\
&+\langle A\xi,W\rangle\langle AY,Z\rangle-\langle A\xi,Z\rangle\langle
AY,W\rangle=0, 
\end{eqnarray*}
for any $Y\perp \xi $, $Z\in \mathcal{H}_{1}(1)$ and $W\in \mathcal{H}%
_{1}(-1)$. Using Lemma \ref{lem:theta}(d) and Lemma \ref{lem:theta_2}(b), we
obtain 
\begin{equation}
\mu \langle V,W\rangle \langle AY,Z\rangle -\beta \langle U,Z\rangle \langle
AY,W\rangle -2\sum_{a=2}^{3}\langle \phi \xi _{a},Y\rangle \langle \phi
_{a}Z,W\rangle =0,  \label{ap2}
\end{equation}%
for any $Y\perp \xi $, $Z\in \mathcal{H}_{1}(1)$ and $W\in \mathcal{H}%
_{1}(-1)$. If $Z\perp U$ and $W\perp V$, then $\sum_{a=2}^{3}\langle \phi
\xi _{a},Y\rangle \langle \phi _{a}Z,W\rangle =0$, for any $Y\perp\xi$.
Since $\phi \xi _{2}$ and $\phi\xi _{3}$ are linearly independent, 
\begin{equation}
\langle \phi _{a}Z,W\rangle =0,\quad a\in \{2,3\},  \label{jati}
\end{equation}
for any $Z\in \mathcal{H}_{1}(1)$ $(\perp U)$ and $W\in \mathcal{H}_{1}(-1)$ 
$(\perp V)$.

If $\dim \mathcal{H}_{1}(1)\geq 4$ the above equation implies that $\phi
_{a} $ is not a monomorphism on $\mathcal{H}_{1}(1)$, which contradicts
Lemma \ref{lem:theta_2}(d). Hence, we conclude that $\dim \mathcal{H}%
_{1}(1)=\dim \mathcal{H}_{1}(-1)=2$, and $\mathcal{H}_{1}(1)=\mathbb{C}U$
and $\mathcal{H}_{1}(-1)=\mathbb{C}V$.

The equation (\ref{jati}) directly implies $\langle \phi _{a}\phi U,\phi V\rangle=0$, for $a\in \{2,3\}$. Hence, by Lemma \ref{lem:theta_2}(d), $\phi_a\phi U=\pm V$, for $a\in \{2,3\}$.
We can express 
\be \phi_a\phi U=\epsilon_a V,\quad \epsilon_a\in\{1,-1\}.\label{eqn:80}\ee


Next, by putting $Z=\phi U$ and $W=V$ in (\ref{ap2}), and using (\ref{eqn:80}) we obtain 
\begin{equation*}
\mu A\phi U=2\sum_{a=2}^{3}\epsilon_a\phi \xi _{a},
\end{equation*}
which implies that $A\phi U\perp \mathcal{H}$.

Putting $Y=U,Z=\phi U, X=V$ in (\ref{ap}) and using Lemma \ref{lem:theta_2}%
(b), we get 
\begin{equation}
R(V,U)\phi U=0.  \notag
\end{equation}
Using Lemma \ref{lem:theta_2}(b), Lemma \ref{lem:theta_2}(d) and the Gauss
equation (\ref{ap1}), we deduce from the above equation that 
\begin{equation}
\phi V-\phi_1V+\sum_{a=2}^{3}\langle V,\phi_a\phi U\rangle\phi_a U=0.
\label{ap3}
\end{equation}
Since $V\in \mathcal{H}_1(-1)$, we have $\theta_1V=-V$, which implies $\phi
V=\phi_1V$. Hence, from (\ref{ap3}) we conclude that $\sum_{a=2}^{3}\langle
V,\phi_a\phi U\rangle\phi_aU=0$.
This contradicts  (\ref{eqn:80}) and the orthogonolity of $\phi_2U$ and $\phi_3U$. Accordingly, this
case cannot occur.

\textbf{Case 2: $\xi\in\mathfrak{D}$ on $M$.}

In this case, we have each $\eta (\xi _{a})=0$ for all $a$, everywhere. It
follows from (\ref{eqn:jacobi}) that the normal Jacobi operator has three
constant eigenvalues $0$, $4$ and $1$ at each point of $M$ with eigenspaces 
\begin{equation*}
T_{0}=\text{Span}\{\phi \xi _{a}:~a=1,2,3\},\quad T_{4}=\text{Span}\{\xi
,\xi _{a}:~a=1,2,3\},\quad T_{1}=\mathcal{H}
\end{equation*}%
respectively.

If we put $X$, $Y\perp \xi $ and $Z=\xi $ in (\ref{ap5}) then $\hat{R}%
_{N}R(X,Y)\xi =4R(X,Y)\xi $ and so $R(X,Y)\xi \in T_{4}\ominus \mathbb{R}\xi
=\text{Span}\{\xi _{a}:~a=1,2,3\}$, for any $X,Y\perp \xi $. Hence, from
Lemma \ref{lem:dga}, it follows that $\langle R(X,Y)\xi ,\phi \xi
_{a}\rangle =0$, for any $X,Y\perp \xi $. Furthermore, using Lemma \ref%
{lem:theta}(d) and the Gauss equation (\ref{ap1}), we obtain 
\begin{equation*}
\eta (AY)\langle AX,\phi \xi_{a}\rangle -\eta (AX)\langle AY,\phi \xi
_{a}\rangle -2\langle \phi _{a}X,Y\rangle =0,
\end{equation*}
for any $X,Y\in \mathcal{H}$ and $a\in \{1,2,3\}$. This equation, together
with Lemma~\ref{lem:Xua}, yields $\langle \phi _{a}X,Y\rangle =0$, for any $%
X,Y\in \mathcal{H}$. This is a contradiction and the proof is completed.


\section{Real hypersurfaces with recurrent $\hat R_N$ in $G_2(\mathbb{C}%
_{m+2})$}

In this section, we show that there does not exist any real hypersurface
with recurrent normal Jacobi operator in $G_{2}(\mathbb{C}_{m+2})$. We first
discuss the ideas of recurrence and semi-parallelism in a general setting. 

Let $M$ be a Riemannian manifold and $\mathcal{E}_j$ a Riemannian vector
bundle over $M$ with linear connection $\nabla^j$, $j\in\{1,2\}$. It is
known that $\mathcal{E}^*_1\otimes \mathcal{E}_2$ is isomorphic to the
vector bundle $Hom(\mathcal{E}_1,\mathcal{E}_2)$, consisting of
homomorphisms from $\mathcal{E}_1$ into $\mathcal{E}_2$. We denote by the
same $\langle,\rangle$ the fiber metrics on $\mathcal{E}_1$ and $\mathcal{E}%
_2$ as well as that induced on $Hom(\mathcal{E}_1,\mathcal{E}_2)$. The
connections $\nabla^1$ and $\nabla^2$ induce on $Hom(\mathcal{E}_1,\mathcal{E%
}_2)$ a connection $\bar \nabla$, given by 
\begin{equation*}
(\bar\nabla_XF)V=(\bar\nabla F)(V;X)=\nabla^2_XFV-F\nabla^1_XV 
\end{equation*}
for any vector field $X$ tangent to $M$, cross sections $V$ in $\mathcal{E}_1
$ and $F$ in Hom$(\mathcal{E}_1,\mathcal{E}_2)$.

A section $F$ in $Hom(\mathcal{E}_1,\mathcal{E}_2)$ is said to be \emph{%
recurrent} if there exists 1-form $\tau$ such that $\bar\nabla F=F\otimes
\tau$. We may regard parallelism as a special case of recurrence, that is,
the case $\tau=0$. Let $\bar R$, $R^1$ and $R^2$ be the curvature tensor
corresponding to $\bar\nabla$, $\nabla^1$ and $\nabla^2$ respectively. Then
we have 
\begin{equation*}
(\bar R\cdot F)(V;X,Y)=(\bar R(X,Y)F)V=R^2(X,Y)FV-FR^1(X,Y)V 
\end{equation*}
for any $X,Y\in TM$, $V\in\mathcal{E}_1$ and $F\in Hom(\mathcal{E}_1,%
\mathcal{E}_2)$.

We have the following result from \cite{nonflat}:

\begin{lem}
\cite{nonflat}\label{nonflat1} Let $M$ be a Riemannian manifold, $\mathcal{E}%
_j$ a Riemannian vector bundle over $M$, $j\in\{1,2\}$ and $F$ a section in $%
Hom(\mathcal{E}_1,\mathcal{E}_2)$. If $F$ is recurrent then $F$ is
semi-parallel.
\end{lem}

From Lemma \ref{nonflat1} and Corollary \ref{cor1.5} we obtain the following:

\begin{cor}\label{4.1}
There does not exist any real hypersurface with recurrent normal Jacobi
operator in $G_{2}(\mathbb{C}_{m+2})$, $m\geq3$.
\end{cor}

As a corollary we have the following:

\begin{cor}
There does not exist any real hypersurface with parallel normal Jacobi
operator in $G_{2}(\mathbb{C}_{m+2})$, $m\geq3$.
\end{cor}

\end{document}